\documentclass{article}

\newtheorem{fed}{\textbf{Definition}}[section]
\newtheorem{thm}[fed]{\textbf{Theorem}}

\usepackage{amssymb,bbm,graphicx,epsfig,psfrag,epic,eepic,latexsym}
\usepackage{amsmath}
\usepackage{mathrsfs}
\usepackage[dvips]{color}

\begin{document}
\title{Periodic orbits in time-dependent planar Stark-Zeeman systems}
\author{Urs Frauenfelder}
\maketitle

\begin{abstract}
Time-dependent Stark-Zeeman systems describe the motion of an electron attracted by a proton subject to a magnetic and a time-dependent electric field.
For instance the study of the dynamics of a gateway around the moon which is subject to the joint attraction of the moon, the earth and the sun leads to 
time-dependent Stark-Zeeman systems. In the time-dependent case there is no preserved energy. Therefore collisions cannot be regularized by blowing up the 
energy hypersurface. A new regularization technique of blowing up instead of the energy hypersurface the loop space was recently discovered by Barutello, Ortega, and Verzini. In this article we explain how this new regularization technique can be applied to the study of periodic orbits in time-dependent planar Stark-Zeeman systems. Since the regularization by blowing-up the loop space is non-local the regularized periodic orbits will not satisfy an ODE anymore but a delay equation. 
\end{abstract}

\section{Introduction}

One of the great challenges for humanity is to settle on the moon. For that purpose it would be useful to have a gateway to the moon, namely a space station on a periodic orbit around the moon. 
The dynamics around the moon is highly complex. If one just takes into account the gateway and the moon one merely has to solve a Newtonian two-body problem whose solution is known since Kepler. If one takes into account the earth as well a good model for the gateway is the restricted three-body problem. By considering this problem in a rotating frame where the earth and the moon are at rest the Hamiltonian for the restricted three-body problem is autonomous, i.e., independent of time. In particular, one has preservation of energy, namely the Hamiltonian is preserved under its one Hamiltonian flow. Recently a lot of progress has been made in understanding the global structure of the network of periodic orbits in the restricted three-body problem
\cite{aydin-bathkin}. To understand how this network is connected one has to take into account as well collisions. Although the actual gateway should have no collisions with the moon, families do in general not end in collisions but can be continuously extended over collisions by regularizing them. To regularize two-body collisions in the autonomous case one can blow-up the energy surface. This is a classical topic in celestial mechanics. How these regularizations are related to contact topology is explained in \cite{frauenfelder-vankoert}.
\\ \\
In a neighbourhood of the moon one does not just feel the gravitational force of the earth but as well on a similar scale the gravitational force of the sun. This means that one should consider the gateway as a restricted four-body problem where the gateway is attracted by three primaries, the moon, the earth, and the sun. A model for this scenario is the bicircular model which is for instance described in \cite{koon-lo-marsden-ross}. Different from the restricted three body problem which starts with an actual solution of the Newtonian two-body problem the bicircular four-body problem uses only an approximative solution to the Newtonian three-body problem where the sun and the center of mass of the earth-moon system rotate around there common center of mass in circles while the moon and the earth rotate additionally around their common center of mass in a circle. In further contrast to the restricted three-body problem one cannot achieve by the transition to rotating coordinates that all three primaries are at rest. However, one can succeed that the earth and the moon are at rest and the sun is than orbiting periodically with period one month. Hence the Hamiltonian for the gateway is not autonomous but depends periodically on time. 
Therefore there is no preserved energy anymore and in the search for periodic orbits one has to go for periodic orbits of period an integer multiple of a month. 
\\ \\
Periodic orbits of fixed period are generically isolated. However, to actually find them and relate them to the global network \cite{aydin-bathkin} one should consider a homotopy which switches on the sun. The additional homotopy parameter has than the effect that the periodic orbits again appear in families. But members of such families might collide with the moon. 
\\ \\
In contrast to the autonomous case in the time-dependent case one does not have energy hypersurfaces anymore since there is no preserved energy. Hence to regularize the system one cannot blow-up the energy hypersurface. A new regularization technique for non-autonomous Hamiltonian systems in the context of forced Kepler problems was recently discovered by Barutello, Ortega, and Verzini in \cite{barutello-ortega-verzini}. In contrast to blowing-up the energy hypersurface one blows-up the loop space. This nonlocal regularization technique was developed further in \cite{cieliebak-frauenfelder-volkov, frauenfelder-weber}.
\\ \\
The goal of this note is to show how the new regularization technique by blowing up the loop space can be applied to the bicircular restricted four-body problems and problems of a similar structure like it. For that purpose we explain in Section~\ref{stark} the notion of a time-dependent Stark-Zeeman system. In the time-independent case this notion was introduced in \cite{cieliebak-frauenfelder-vankoert}. It considers the Coulomb problem or equivalently the Kepler problem subject to an additional electric and magnetic force. In the time-dependent case we allow in addition that the electric force can depend periodically on time. The bicircular restricted four-body problem is an instance of a time-dependent Stark-Zeeman system. Since one considers the bicircular restricted four-body problem in a rotating frame the system will be subject to additional forces like Coriolis force and centrifugal force. In contrast to the centrifugal force and the gravitational force the Coriolis force is velocity dependent and therefore can be modelled analogously as the Lorentz force with the help of a magnetic field. 
The electric force is a combination of the centrifugal force and the gravitational forces of the earth and the sun. In particular the electric force depends periodically on time since the sun is on periodic motion. 
\\ \\
In Section~\ref{blow} we explain how the new regularitation technique of blowing up the loop space which was first applied in \cite{barutello-ortega-verzini} to forced Kepler problems translates to periodic orbits in time-dependent Stark-Zeeman systems. As the classical Levi-Civita regularization \cite{levicivita} one uses in this regularization technique the complex squaring map. In the classical Levi-Civita regularization one changes time by multiplying the energy hypersurface which a positive function. Since in non-autonomous systems there is no preserved energy this technique does not work. Instead of that one reparametrizes each loop individually. In particular, this regularization technique is not local but depends on the knowledge of the whole loop. 
\\ \\
Since the regularization depends on the whole loop the resulting equation is not local anymore and a second order delay equation. In Section~\ref{critic} we derive this equation. 
\\ \\
In Section~\ref{regular} we discuss the relation between the solutions of the unregularized original second order ODE and the regularized second order delay equation for noncollisional solutions. If the solution has no collisions already in the unregularized case there is a variational approach to detect them. The regularized solution is than just obtained as critical point of the pulled-back action functional. Hence it is clear that a regularized noncollisional solution corresponds to a nonregularized one. Nevertheless it is interesting to study what happens precisely when one transforms the regularized second order delay equation back. By doing that we will discover how even the unregularized solutions arise as solutions of a delay equation.
\\ \\
In Section~\ref{coll} we explain how in the case of collisions solutions of the regularized second order delay equation correspond to collisional solutions of the original second order ODE. 
The regularized solutions in the case of collisions are still smooth and collisions just appear when the solution goes through the origin. On the other hand solutions for the original second order ODE in the case of collisions are merely continuous and its derivative explodes at collisions. There is no variational approach anymore for the unregularized problem. However, the regularized solutions still appear as critical points of an action functional. On the other hand 
we see in this Section that after transforming back we still obtain a similar delay equation for the unregularized solutions as we obtained in Section~\ref{regular}. We than analysis how the solutions of this delay equation corresponds to collisional solutions of the original second order ODE.
\\ \\
In Section~\ref{correspondence} we discuss the precise correspondence between the regularized solutions and the original collisional solutions. The main findings of this note we summarize in Theorem~\ref{main}. The reason why the correspondence is not completely direct is that the complex squaring map is a two-to-one cover. On the other hand this also has an interesting consequence. Although we cannot make sense anymore of the winding number for collisional solutions of the original second order ODE we still can make sense of the winding number modulo two. 
\\ \\
A different approach for regularizing collisions in the case of non-autonomous Hamiltonians depending periodically on time would be the extended phase space approach. In the case of forced Kepler problems this was carried out by Zhao \cite{zhao}. In fact it might be advantageous to have several approaches to regularizations available. The ultimate goal is to use them to homotop 
periodic solutions from the restricted three-body problem to the bicircular restricted four-body problem by switching on the sun. This homotopy might depend on how the sun is switched on so that there is the possibility that some monodromy is discovered by looking at different homotopies. 
\\ \\
There are various further directions this work can be extended. Although for the gateway one is mainly interested for periodic orbits close to the moon and therefore the major concern are collisions with the moon in principle it could happen that there are families which have members which collide with the earth. To take collisions with the earth into account as well one could look at time-dependent, two-center Stark-Zeeman systems. In the autonomous case two-center Stark-Zeeman systems were introduced in \cite{cieliebak-frauenfelder-zhao}. There it is as well explained how with the help of a description given in \cite{waldvogel} the complex squaring map can be modified to give rise to Birkhoff regularization which regularizes both primaries simultaneously. In the time-dependent case there should also exist a nonlocal version of Birkhoff's regularization. In \cite{waldvogel} it is as well explained how quaternions instead of complex numbers can be applied to regularize collisions in the autonomous spatial case which gives rise to the Kustaanheimo-Stiefel regularization extending the Levi-Civita regularization. There should also exist a nonlocal version of the Kuustaanheimo-Stiefel regularization for spatial time-dependent Stark-Zeeman systems. In the case of forced Kepler problems such a regularization was already studied in \cite{barutello-ortega-verzini}. A further direction is to take into account fidelity models than the bicircular restricted three-body problems as Scheeres version of the Hill four-body problem \cite{scheeres} or the quasi-bicircular problem \cite{andreu}. These models take additionally into account the varying distance between earth and moon and hence lead to pulsating coordinates.
\\ \\  
\emph{Acknowledgements:} The author acknowledges partial support by DFG project FR 2637/6-1.

\section{Time-dependent Stark-Zeeman systems}\label{stark}
We assume that $\mathfrak{Q}_0$ is an open subset of $\mathbb{C}$ containing the origin and 
$$\mathfrak{Q}=\mathfrak{Q}_0 \setminus \{0\}$$
is the open subset of $\mathbb{C}$ obtained by removing the origin from $\mathfrak{Q}_0$. Abbreviate by
$S^1=\mathbb{R}/\mathbb{Z}$ the circle. Suppose that
$$E \colon \mathfrak{Q}_0 \times S^1 \to \mathbb{R}$$
is a smooth function. For $t \in S^1$ we abbreviate 
$$E_t \colon \mathfrak{Q}_0 \to \mathbb{R}, \quad q \mapsto E(q,t)$$
and refer to it as the \emph{electric potential} depending periodically on time. We add to the electric potential the Coulomb potential to get for $t \in S^1$ the time-dependent potential
$$V_t \colon \mathfrak{Q} \to \mathbb{R}, \quad
q \mapsto E_t(q)-\frac{1}{|q|}.$$
We abbreviate for the Hamiltonian
$$H_t \colon T^* \mathfrak{Q}, \colon (q,p) \mapsto
\frac{1}{2}|p|^2+V_t(q),$$
which consists of kinetic and potential energy.
\\ \\
To incorporate a magnetic field to the Hamiltonian system we can twist the standard symplectic form on the cotangent bundle. For that purpose,  assume that 
$$\sigma=Bdq_1 \wedge dq_2 \in \Omega^2(\mathfrak{Q}_0)$$
is a two-form on $\mathfrak{Q}_0$. 
We refer to
$$B \in C^\infty(\mathfrak{Q}_0,\mathbb{R})$$
as the \emph{magnetic field}.
Since $\mathfrak{Q}_0$ is two-dimensional $\sigma$ is necessarily closed, and  because the second de Rham cohomology of an open subset in $\mathbb{C}$ is trivial, it follows that $\sigma$ is even exact. This means that there exists a one-form $A \in \Omega^1(\mathfrak{Q}_0)$ such that 
$$\sigma=dA.$$ 
If we write $A=A_1dq_1+A_2dq_2$, then we have
$$B=\mathrm{rot}A=\frac{\partial A_2}{\partial q_1}-\frac{\partial A_1}{\partial q_2}.$$
The Liouville one-form $\lambda \in \Omega^1(T^* \mathbb{C})$ at a point
$(q,p)=(q_1+iq_2,p_1+ip_2) \in \mathbb{C} \oplus \mathbb{C}=T^* \mathbb{C}$ is given by
$$\lambda=p_1 dq_1+p_2dq_2.$$
Its exterior differential is the standard symplectic form
on $T^*\mathbb{C}$ given by
$$\omega=d\lambda=dp_1 \wedge dq_1+dp_2 \wedge dq_2.$$
Denote by $\pi \colon T^* \mathbb{C} \to \mathbb{C}, \quad (q,p) \mapsto q$ the footpoint projection. Using $A$ we can twist the Liouville one form to the 
one-form
$$\lambda_A:=\lambda +\pi^*A \in \Omega^1(T^* \mathfrak{Q}_0).$$
Its exterior differential is than the twisted symplectic form
$$d\lambda_A=d\lambda+\pi^*dA=\omega+\pi^*\sigma=:\omega_\sigma \in \Omega^2(T^* \mathfrak{Q}_0).$$
We define the twisted Hamiltonian vector field $X_{H_t}^\sigma$ on $T^* \mathfrak{Q}$ of the Hamiltonian $H_t$ implicitly by the requirement that
$$dH_t=\omega_\sigma(\cdot, X_{H_t}^\sigma).$$
A periodic orbit $x \in C^\infty(S^1,T^*\mathfrak{Q})$ is a solution of the first order ODE
\begin{equation}\label{ham}
\dot{x}(t)=X_{H_t}^\sigma(x(t)), \quad t \in S^1.
\end{equation}
Writing $x=(q,p)$ we can rewrite this first order ODE for $x$ as a second order ODE for $q$. Using the formulas
$$dH_t=p_1dp_1+p_2dp_2+\frac{\partial V_t}{\partial q_1}dq_1+\frac{\partial V_t}{\partial q_2}dq_2, \quad \omega_\sigma=dp_1\wedge dq_1+dp_2\wedge dq_2+Bdq_1\wedge dq_2$$
we can write the Hamiltonian vector field explicitly as
$$X_{H_t}^\sigma=p_1\frac{\partial}{\partial q_1}
+p_2\frac{\partial}{\partial q_2}+\bigg(Bp_2-\frac{\partial V_t}{\partial q_1}
\bigg)\frac{\partial}{\partial p_1}-
\bigg(Bp_1+\frac{\partial V_t}{\partial q_2}\bigg)
\frac{\partial}{\partial p_2}.$$
Therefore (\ref{ham}) is equivalent to the following system of equations
\begin{eqnarray}\label{ode1}
\dot{q}_1&=&p_1\\ \nonumber
\dot{q}_2&=&p_2\\ \nonumber
\dot{p}_1&=&Bp_2-\frac{\partial_t V}{\partial q_1}\\ \nonumber
\dot{p}_2&=&-Bp_1-\frac{\partial V_t}{\partial q_2}.	
\end{eqnarray}	
Using for the gradient
$$\nabla V_t(q)=\frac{q}{|q|^3}+\nabla E_t(q)$$
and taking into account again complex notation we can rewrite the above system (\ref{ode1}) of a first order ODE as the second order ODE for 
$q \in C^\infty(S^1,\mathfrak{Q})$
\begin{equation}\label{ode2}
\ddot{q}=-B(q)i\dot{q}-\frac{q}{|q|^3}-\nabla E_t(q).
\end{equation}
We abbreviate the free loop space of $\mathfrak{Q}$ by
$$\mathfrak{LQ}=C^\infty(S^1,\mathfrak{Q}).$$
This is an open subset of the free loop space $\mathfrak{L}\mathbb{C}$ of $\mathbb{C}$. We abbreviate
$$\langle \cdot,\cdot \rangle \colon \mathfrak{L}\mathbb{C} \times \mathfrak{L}\mathbb{C} \to \mathbb{R}$$
the $L^2$-inner product given for $\xi, \eta \in \mathfrak{L}\mathbb{C}$ by
$$\langle \xi,\eta \rangle=\int_0^1 \langle \xi(t),\eta(t)\rangle dt=\int_0^1 \mathrm{Re}\big(i\xi(t) \bar{\eta}(t)\big)dt$$
where $\bar{\eta}$ denotes the complex conjugate. By
$$|| \cdot|| \colon \mathfrak{L}\mathbb{C} \to [0,\infty), \quad \xi \mapsto \sqrt{\langle \xi,\xi\rangle}$$
we abbreviate the $L^2$-norm. Using this notation periodic orbits satisfying (\ref{ode2}) can be described variationally using the Lagrangian formalism by the action functional
$$\mathfrak{F} \colon \mathfrak{LQ}\to \mathbb{R}, \quad q \mapsto \frac{1}{2}||\dot{q}||^2-\int_{S^1}q^*A-\int_0^1 V_t(q(t))dt.$$
Indeed, suppose that $q$ is a critical point of $\mathfrak{F}$ and $\xi \in T_q \mathfrak{LQ}=\mathfrak{L}\mathbb{C}$ is a tangent vector. Than
\begin{eqnarray*}
0&=&d\mathfrak{F}(q)\xi\\
&=&\langle \dot{q},\dot{\xi}\rangle-\int_0^1 dA(q(t))\big(\dot{q}(t),\xi(t)\big)dt-\int_0^1 dV_t(q(t))\xi(t)dt\\
&=&-\langle \ddot{q},\xi\rangle-\langle Bi\dot{q},\xi\rangle-\langle \nabla V,\xi\rangle\\
&=&-\langle \ddot{q}+Bi\dot{q}+\nabla V,\xi\rangle,
\end{eqnarray*}
so that $q$ is a solution of (\ref{ode2}).

\section{Blowing up the loop space}\label{blow}

We abbreviate by
$$\varsigma \colon \mathbb{C} \to \mathbb{C}, \quad z \mapsto z^2$$
the complex squaring map. By
$$\mathfrak{Z}:=\varsigma^{-1}(\mathfrak{Q}) \subset \mathbb{C}$$
we denote the preimage of $\mathfrak{Q}$ under the complex squaring map. Then $\mathfrak{Z}$ is an open subset of $\mathbb{C}$ which does not contain the origin, but whose closure contains the origin. For later purpose we also introduce 
$$\mathfrak{Z}_0:=\varsigma^{-1}(\mathfrak{Q}_0)=\mathfrak{Z} \cup \{0\} \subset \mathbb{C}.$$
For $z \in \mathfrak{LZ}$ we introduce the map
$$t_z \colon S^1 \to S^1, \quad \tau \mapsto \frac{\int_0^\tau |z(s)|^2ds}{||z||^2}.$$
Its derivative
$$t'_z(\tau)=\frac{|z(\tau)|^2}{||z||^2}$$
is positive, since $z(\tau)\neq 0$ for every $\tau \in S^1$. In particular, $t_z$ is a diffeomorphism and we denote by
$$\tau_z=t_z^{-1} \colon S^1 \to S^1$$
its inverse. Its derivative at $t \in S^1$ is given by
\begin{equation}\label{dertau}
\dot{\tau}_z(t)=\frac{1}{t'_z(\tau_z(t))}=\frac{||z||^2}{|z(\tau_z(t))|^2}.
\end{equation}
We further define
$$\Sigma \colon \mathfrak{LZ} \to \mathfrak{LQ}, \quad z \mapsto q_z$$
where 
$$q_z(t):=z^2(\tau_z(t)), \quad t \in S^1.$$
We next compute
$$\Sigma^* \mathfrak{F}=\mathfrak{F} \circ \Sigma \colon \mathfrak{LZ} \to \mathbb{R}.$$
In order to do that we first compute
\begin{equation}\label{der1}
\dot{q}_z(t)=2z(\tau_z(t))z'(\tau_z(t))\dot{\tau_z}(t)=2||z||^2\frac{z'(\tau_z(t))}{\bar{z}(\tau_z(t))}
\end{equation}
to obtain
\begin{eqnarray}\label{kin}
||\dot{q}_z||^2&=&4\int_0^1 ||z||^4\frac{|z'(\tau_z(t))|^2}{|z(\tau_z(t))|^2}dt\\ \nonumber
&=&4\int_0^1||z||^4\frac{|z'(\tau)|^2}{|z(\tau)|^2}t'_z(\tau)d\tau\\ \nonumber
&=&4||z||^2||z'||^2.
\end{eqnarray}
We further note
\begin{equation}\label{pot}
\int_0^1 \frac{1}{|q_z(t)|}dt=\int_0^1 \frac{1}{|z(\tau)|^2}t'_z(\tau)d\tau=\frac{1}{||z||^2}.
\end{equation}
and 
\begin{eqnarray}\label{ener}
\int_0^1 E_t(q_z(t))dt&=&\frac{1}{||z||^2}\int_0^1 E_{t_z(\tau)}\big(z^2(\tau)\big)|z(\tau)|^2d\tau.
\end{eqnarray}
Since integrating a one-form does not depend on the parametrization we obtain for the pullback of $\mathfrak{F}$ at $z \in \mathfrak{LZ}$
\begin{eqnarray}\label{pullback}
& &\Sigma^*\mathfrak{F}(z)=\mathfrak{F}(q_z)\\ \nonumber
&=&2||z||^2||z'||^2-\int_{S^1} z^*\varsigma^* A-\frac{1}{||z||^2}-\frac{1}{||z||^2}\int_0^1 E_{t_z(\tau)}\big(z^2(\tau)\big)|z(\tau)|^2d\tau.
\end{eqnarray}
We observe that this functional extends to the following open subset of the free loop space of $\mathfrak{LZ}_0$
$$\mathfrak{L^*Z}_0:=\big\{z \in \mathfrak{LZ}_0: ||z|| \neq 0\big\}.$$
Indeed define
$$\mathcal{F} \colon \mathfrak{L^*Z}_0 \to \mathbb{R}$$
for $z \in \mathfrak{L^*Z}_0$ by the formula above, namely
$$\mathcal{F}(z)=2||z||^2||z'||^2-\int_{S^1} z^*\varsigma^* A-\frac{1}{||z||^2}-\frac{1}{||z||^2}\int_0^1 E_{t_z(\tau)}\big(z^2(\tau)\big)|z(\tau)|^2d\tau$$
so that we have
$$\mathcal{F}|_{\mathfrak{LZ}}=\Sigma^*\mathfrak{F}.$$

\section{Critical points}\label{critic}

In this section we determine the critical point equation for the functional $\mathcal{F}$. Suppose that $ z \in \mathfrak{L^*Z}_0$ is a critical point of
$\mathcal{F}$. This means that for every $\zeta \in T_z \mathfrak{L^*Z}=\mathfrak{L}\mathbb{C}$ we have
$$d\mathcal{F}(z)\zeta=0.$$
In order to derive an explicit description of the differential of the functional $\mathcal{F}$ we first write it as the sum of four terms. We first set
$$\mathcal{K} \colon \mathfrak{L^*Z}_0 \to \mathbb{R}, \quad z \mapsto 2||z||^2||z'||^2$$
and refer to this functional as the \emph{kinetic part}. The functional
$$\mathcal{A} \colon \mathfrak{L^*Z}_0 \to \mathbb{R}, \quad z \mapsto \int_{S^1} z^*\varsigma^* A$$
we refer to as the \emph{magnetic part}, the functional
$$\mathcal{C} \colon \mathfrak{L^*Z}_0 \to \mathbb{R}, \quad z \mapsto \frac{1}{||z||^2}$$
as the \emph{Coulomb part}, and the functional
$$\mathcal{E} \colon \mathfrak{L^*Z}_0 \to \mathbb{R}, \quad z \mapsto \frac{1}{||z||^2}\int_0^1 E_{t_z(\tau)}\big(z^2(\tau)\big)|z(\tau)|^2d\tau$$
as the \emph{electric part}. Hence we can write
$$\mathcal{F}=\mathcal{K}-\mathcal{A}+\mathcal{C}-\mathcal{E}.$$
For $z \in \mathfrak{L^*Z}_0$ and $\zeta \in \mathfrak{L}\mathbb{C}$ the differential of the kinetic part computes to be
\begin{equation}\label{dk}
d\mathcal{K}(z)\zeta=4||z'||^2\langle z,\zeta\rangle +4||z||^2\langle z',\zeta'\rangle=4||z'||^2\langle z,\zeta\rangle -4||z||^2\langle z'',\zeta\rangle.
\end{equation}
Before deriving the potential of the magnetic part, we first compute the pull-back of the two-form $\sigma=Bdq_1 \wedge dq_2$ under the squaring map $\varsigma$. In order to do that we write
$$q_1+iq_2=q=z^2=(z_1+iz_2)^2=z_1^2-z_2^2+2iz_1z_2$$
so that we obtain
$$q_1=z_1^2-z_2^2, \quad q_2=2z_1z_2$$
and therefore the volume form $dq_1 \wedge dq_2$ pulls-back under $\varsigma$ to 
$$dq_1 \wedge dq_2=4(z_1dz_1-z_2dz_2)\wedge(z_2dz_1+z_1dz_2)=4(z_1^2+z_2^2)dz_1 \wedge dz_2=4|z|^2dz_1 \wedge dz_2.$$
Hence we have
$$\varsigma^*\sigma(z)=4|z|^2 B(z^2)dz_1\wedge dz_2.$$
Hence we get for the differential of the magnetic part
\begin{equation}\label{da}
d\mathcal{A}(z)\zeta=4\big\langle|z|^2B(z^2)iz',\zeta\big\rangle.
\end{equation}
The differential of the Coulomb part is
\begin{equation}\label{dc}
d\mathcal{C}(z)\zeta=-\frac{2\langle z,\zeta\rangle}{||z||^4}.
\end{equation}
To compute the differential of the electric part we introduce the following auxiliary notions.
We define
$$\mathcal{E}^1 \colon \mathfrak{L^*Z}_0 \to \mathbb{R}$$
for $z \in \mathfrak{L^*Z}_0$ by
$$\mathcal{E}^1(z)=\frac{1}{||z||^4}\int_0^1 \bigg(\int_0^\sigma |z(s)|^2ds\bigg) \dot{E}_{t_z(\sigma)}(z^2(\sigma))|z(\sigma)|^2 d\sigma.$$
We write  the gradient as a complex number
$$\nabla E(x+iy):=\frac{\partial E}{\partial x}(x+iy)+i\frac{\partial E}{\partial y}(x+iy)
\in \mathbb{C}$$
and for its complex conjugate
$$\bar{\nabla} E(x+iy):=\frac{\partial E}{\partial x}(x+iy)-i\frac{\partial E}{\partial y}(x+iy)
\in \mathbb{C}.$$
Using this notation we introduce the following three vector field
$$\varepsilon_k \colon \mathfrak{L^*Z}_0 \to \mathfrak{L}\mathbb{C}, \quad k \in \{1,2,3\}.$$
Namely for $z \in \mathfrak{L^*Z}_0$ and $\tau \in S^1$ we set
\begin{eqnarray*}
\varepsilon_1(z)(\tau)&=&\frac{1}{||z||^2}\bigg(\int_\tau^1\dot{E}_{t_z(\sigma)}(z^2(\sigma))
|z(\sigma)|^2d\sigma\bigg)z(\tau)\\
\varepsilon_2(z)(\tau)&=&|z(\tau)|^2\nabla E_{t_z(\tau)}(z^2(\tau))\bar{z}(\tau)\\
\varepsilon_3(z)(\tau)&=&E_{t_z(\tau)}(z^2(\tau))z(\tau).
\end{eqnarray*} 
Adding these three vector fields together we get the vector field
$$\varepsilon=\varepsilon_1+\varepsilon_2+\varepsilon_3 \colon \mathfrak{L^*Z}_0 \to \mathfrak{L}\mathbb{C}.$$
Taking into account the formula
\begin{eqnarray*}
	(dt_z \zeta)(\tau)&=&\frac{2\int_0^\tau \mathrm{Re}\big(\bar{z}(\sigma)\zeta(\sigma)\big)d\sigma}{||z||^2}-
	\frac{2\langle z,\zeta\rangle \int_0^\tau |z(\sigma)|^2d\sigma}{||z||^4}
\end{eqnarray*}
we compute
\begin{eqnarray}\label{de}
	d\mathcal{E}(z)\zeta&=&-\frac{2\langle z,\zeta\rangle}{||z||^4}\int_0^1 E_{t_z(\tau)}\big(z^2(\tau)\big)|z(\tau)|^2d\tau\\ \nonumber
	& &+\frac{1}{||z||^2}\int_0^1 \dot{E}_{t_z(\tau)}(z^2(\tau))(dt_z \zeta)(\tau)|z(\tau)|^2 d\tau\\ \nonumber
	& &+\frac{2}{||z||^2}\int_0^1 |z(\tau)|^2\mathrm{Re}\Big(\bar{\nabla} E_{t_z(\tau)}(z^2(\tau))
	z(\tau)\zeta(\tau)\Big)d\tau\\ \nonumber
	& &+\frac{2}{||z||^2}\int_0^1 E_{t_z(\tau)}(z^2(\tau))\mathrm{Re}\big(\bar{z}(\tau)\zeta(\tau)\big)d\tau\\ \nonumber
	&=&-\frac{2\langle z,\zeta\rangle}{||z||^2}\mathcal{E}(z)\\ \nonumber
	& &+\frac{2}{||z||^4}\int_0^1 \dot{E}_{t_z(\tau)}(z^2(\tau))|z(\tau)|^2 \bigg(\int_0^\tau \mathrm{Re}\big(\bar{z}(\sigma)
	\zeta(\sigma)\big)d\sigma \bigg)d\tau\\ \nonumber
	& &-\frac{2\langle z,\zeta\rangle}{||z||^6}\int_0^1 \dot{E}_{t_z(\tau)}(z^2(\tau))|z(\tau)|^2 \bigg(\int_0^\tau |z(\sigma)|^2d\sigma \bigg)d\tau\\ \nonumber
	& &+\frac{2}{||z||^2}\int_0^1 |z(\tau)|^2\mathrm{Re}\Big(\bar{\nabla} E_{t_z(\tau)}(z^2(\tau))
	z(\tau)\zeta(\tau)\Big)d\tau\\ \nonumber
	& &+\frac{2}{||z||^2}\int_0^1 E_{t_z(\tau)}(z^2(\tau))\mathrm{Re}\big(\bar{z}(\tau)\zeta(\tau)\big)d\tau\\ \nonumber
	&=&-\frac{2\langle z,\zeta\rangle}{||z||^2}\mathcal{E}(z)\\ \nonumber
	& &+\frac{2}{||z||^4}\int_0^1\bigg( \int_\sigma^1 \dot{E}_{t_z(\tau)}(z^2(\tau))|z(\tau)|^2 d\tau\bigg)\mathrm{Re}\big(\bar{z}(\sigma)
	\zeta(\sigma)\big)d\sigma\\ \nonumber
	& &-\frac{2\langle z,\zeta\rangle}{||z||^2}\mathcal{E}^1(z)\\ \nonumber
	& &+\frac{2}{||z||^2}\int_0^1 |z(\tau)|^2\mathrm{Re}\Big(\bar{\nabla} E_{t_z(\tau)}(z^2(\tau))
	z(\tau)\zeta(\tau)\Big)d\tau\\ \nonumber
	& &+\frac{2}{||z||^2}\int_0^1 E_{t_z(\tau)}(z^2(\tau))\mathrm{Re}\big(\bar{z}(\tau)\zeta(\tau)\big)d\tau\\
	&=&-\frac{2(\mathcal{E}+\mathcal{E}^1)\langle z,\zeta)}{||z||^2}+\frac{2\langle \varepsilon(z),\zeta\rangle}{||z||^2}.
\end{eqnarray}
A critical point $z$ of the functional $\mathcal{F}$ satisfies
$$d\mathcal{K}(z)\zeta-d\mathcal{A}(z)\zeta+d\mathcal{C}(z)\zeta-d\mathcal{E}(z)\zeta=0,
\quad \forall \zeta \in \mathfrak{L}\mathbb{C}$$
and therefore it follows from (\ref{dk}), (\ref{da}), (\ref{dc}), and (\ref{de}) that $z$ is a solution of the following second order delay equation for every $\tau \in S^1$
\begin{eqnarray}\label{delay}
z''(\tau)&=&-\frac{\varepsilon(z)(\tau)}{2||z||^4}-\frac{|z(\tau)|^2B(z^2(\tau))}{||z||^2}iz'(\tau)\\ \nonumber
& &+\bigg(\frac{||z'||^2}{||z||^2}+\frac{\mathcal{E}(z)+\mathcal{E}^1(z)}{2||z||^4}-\frac{1}{2||z||^6}\bigg)z(\tau).	
\end{eqnarray}

\section{Regular solutions}\label{regular}

If $z \in \mathfrak{LZ} \subset \mathfrak{L^*Z}_0$ is a solution of the critical point equation (\ref{delay}), then $q=\Sigma(z) \in \mathfrak{LQ}$ is a solution of the second order ODE (\ref{ode2}). This is clear, since solutions of (\ref{ode2}) are critical points of the functional $\mathfrak{F} \colon \mathfrak{LQ} \to \mathbb{R}$ and the restriction of the functional $\mathcal{F}$ to $\mathfrak{LZ}$ coincides with the pull-back of $\mathfrak{F}$ under $\Sigma$. Although this is clear it is not really obvious by looking at the equations (\ref{delay}) and (\ref{ode2}). Therefore in this section we check this explicitly. In carrying this out we derive a second order delay equation for $q$ which will help us to understand in the following how solutions of (\ref{delay}) in the complement $\mathfrak{L^*Z}_0 \setminus\mathfrak{LZ}$ correspond to collisional solutions of (\ref{ode2}).
\\ \\
The first derivative of $q=\Sigma(z)=z^2(\tau_z(t))$ was already computed in (\ref{der1})
to be
$$\dot{q}(t)=2||z||^2\frac{z'(\tau_z(t))}{\bar{z}(\tau_z(t))}.$$
Using (\ref{dertau}) we compute for its second derivative
\begin{eqnarray}\label{dotdot}
	\ddot{q}(t)&=&\frac{2||z||^4}{|z(\tau(t))|^2 \bar{z}(\tau(t))}
	\bigg(z''(\tau(t))-\frac{|z'(\tau(t))|^2}{\bar{z}(\tau(t))}\bigg)\\ \nonumber
	&=&\frac{1}{\bar{q}(t)}\bigg(\frac{2||z||^4 z''(\tau_z(t))}{z(\tau_z(t))}-\frac{1}{2}|\dot{q}(t)|^2\bigg).
\end{eqnarray}
To take advantage of (\ref{delay}) we evaluate each term on the righthand side of (\ref{delay}) recalling that $\varepsilon=\varepsilon_1+\varepsilon_2+\varepsilon_3$. We first compute for $\varepsilon_1$
\begin{eqnarray}\label{t1}
\frac{\varepsilon_1(\tau_z(t))}{z(\tau_z(t))}&=&
\frac{1}{||z||^2}\int_{\tau_z(t)}^1\dot{E}_{t_z(\sigma)}(z^2(\sigma))|z(\sigma)|^2d\sigma
\\ \nonumber
&=&\int_t^1 \dot{E}_s(q(s))ds.
\end{eqnarray}	
For $\varepsilon_2$ we derive
\begin{eqnarray}\label{t2}
	\frac{\varepsilon_2(\tau_z(t))}{\bar{q}(t)z(\tau_z(t))}&=&\nabla E_t(q(t)).
\end{eqnarray}	
For $\varepsilon_3$ we get
\begin{equation}\label{t3}
\frac{\varepsilon_3(\tau_z(t))}{z(\tau_z(t))}=E_t(q(t)).	
\end{equation}
We next compute
\begin{eqnarray}\label{t4}
\frac{2||z||^4|z(\tau_z(t))|^2B(q(t))}{\bar{q}(t)z(\tau_z(t))||z||^2}iz'(\tau_z(t))&=&\frac{2||z||^2 B(q(t))}{\bar{z}(\tau_z(t))}iz'(\tau_z(t))\\ \nonumber
&=&B(q(t))i\dot{q}(t).
\end{eqnarray}
If we define the functional
$$\mathfrak{E} \colon \mathfrak{LQ}\to \mathbb{R}, \quad q \mapsto \int_0^1 E_t(q(t))dt$$
then it follows from (\ref{ener}) that
\begin{equation}\label{t5}
\mathcal{E}(z)=\mathfrak{E}(q).
\end{equation}
For $\mathfrak{E}^1(z)$ we compute
\begin{eqnarray*}
\mathcal{E}^1(z)&=&\frac{1}{||z||^4}\int_0^1 \bigg(\int_0^\sigma |z(s)|^2ds\bigg) \dot{E}_{t_z(\sigma)}(z^2(\sigma))|z(\sigma)|^2 d\sigma\\
&=&\int_0^1\bigg(\int_0^s dt\bigg)\dot{E}_s(q(s))ds\\
&=&\int_0^1 s\dot{E}_s(q(s))ds.
\end{eqnarray*}
Hence if we define
$$\mathfrak{E}^1 \colon \mathfrak{LQ}\to \mathbb{R}, \quad q \mapsto \int_0^1 s\dot{E}_s(q(s))ds$$
we can summarize the above computation into the formula
\begin{equation}\label{t6}
	\mathcal{E}^1(z)=\mathfrak{E}^1(q).
\end{equation}
Hence plugging (\ref{delay}) into (\ref{dotdot}) and using (\ref{t1})--(\ref{t6}) together with (\ref{kin}) and (\ref{pot}) we obtain for $q$ the following second order delay equation
\begin{eqnarray}\label{delayq}
& &\ddot{q}(t)+B(q(t))i\dot{q}(t)+\nabla E_t(q(t))\\ \nonumber
&=&\frac{1}{\bar{q}(t)}\bigg(\mathfrak{E}(q)+\mathfrak{E}^1(q)+\frac{||\dot{q}||^2}{2}
-\int_0^1 \frac{1}{|q(s)|}ds\bigg)\\ \nonumber
& &-\frac{1}{\bar{q}(t)}\bigg(
\int_t^1 \dot{E}_s(q(s))ds+E_t(q(t))+\frac{|\dot{q}(t)|^2}{2}\bigg).
\end{eqnarray}
It remains to show that if $q \in \mathfrak{LQ}$ is a solution of the second order delay equation (\ref{delayq}) it is as well a solution of the
second order ODE (\ref{ode2}). For that purpose we set
\begin{equation}\label{bet}
\beta=\frac{\ddot{q}+B(q)i\dot{q}+\nabla E(q)}{q}.
\end{equation}
In view of (\ref{delayq}) we obtain
\begin{eqnarray}\label{beq}
\beta(t) |q(t)|^2&=&
\mathfrak{E}(q)+\mathfrak{E}^1(q)+\frac{||\dot{q}||^2}{2}
-\int_0^1 \frac{1}{|q(s)|}ds\\ \nonumber
& &-\int_t^1 \dot{E}_s(q(s))ds-E_t(q(t))-\frac{|\dot{q}(t)|^2}{2}.
\end{eqnarray}
The righthand side of (\ref{delayq}) is real and therefore
$$\beta(t) \in \mathbb{R}, \quad t \in S^1.$$
Differentiating (\ref{beq}) and using that $\beta$ is real we obtain
\begin{eqnarray*}
2\big\langle q(t),\dot{q}(t)\big\rangle\beta(t)+|q(t)|^2\dot{\beta}(t)&=&\dot{E}_t(q(t))-\dot{E}_t(q(t))\\
& &-\big\langle \nabla E_t(q(t)),\dot{q}(t)\rangle-\big\langle \dot{q}(t),\ddot{q}(t)\big\rangle\\
&=&-\big\langle \nabla E_t(q(t)),\dot{q}(t)\rangle-\big\langle \dot{q}(t),\beta(t) q(t)\big\rangle\\
& &+\big\langle \dot{q}(t),B(q(t))i\dot{q}(t)\big \rangle+\big \langle \dot{q}(t),\nabla E_t(q(t))\big\rangle\\
&=&-\big\langle q(t),\dot{q}(t)\big\rangle.
\end{eqnarray*}
which we rearrange
$$|q|^2\dot{\beta}=-3\langle q,\dot{q}\rangle \beta$$
to conclude
$$\frac{\dot{\beta}}{\beta}=-\frac{3}{2}\frac{\partial_t |q|^2}{|q|^2}.$$
Solving this first order ODE implies that there exists $c \in \mathbb{R}$ such that
$$\ln |\beta|=-\frac{3}{2}\ln |q|^2+c=-\ln |q|^3+c$$
respectively
\begin{equation}\label{beq2}
\beta=\frac{\mu}{|q|^3}
\end{equation}
for $\mu=\pm e^c$.
In view of the definition (\ref{bet}) of $\beta$ we conclude that $q$ is a solution of the following second order ODE 
\begin{equation}\label{mode}
\ddot{q}=\frac{\mu q}{|q|^3}-B(q)i\dot{q}-\nabla E(q).
\end{equation}
It remains to determine $\mu$. For that purpose we first note that combining (\ref{beq}) and (\ref{beq2}) we obtain
\begin{eqnarray*}
\frac{\mu}{|q(t)|}&=&\beta(t)|q(t)|^2\\
&=&\mathfrak{E}(q)+\mathfrak{E}^1(q)+\frac{||\dot{q}||^2}{2}
-\int_0^1 \frac{1}{|q(s)|}ds-\int_t^1 \dot{E}_s(q(s))ds-E_t(q(t))\\
& &-\frac{|\dot{q}(t)|^2}{2}.	
\end{eqnarray*}
Integrating this equation we obtain
\begin{eqnarray*}
\mu\int_0^1 \frac{1}{|q(t)|}dt&=&\mathfrak{E}(q)+\mathfrak{E}^1(q)+\frac{||\dot{q}||^2}{2}
-\int_0^1 \frac{1}{|q(t)|}dt-\int_0^1 \int_t^1\dot{E}_s(q(s))dsdt\\
& &-\int_0^1E_t(q(t))dt-\int_0^1\frac{|\dot{q}(t)|^2}{2}dt\\
&=&\mathfrak{E}(q)+\mathfrak{E}^1(q)+\frac{||\dot{q}||^2}{2}
-\int_0^1 \frac{1}{|q(t)|}dt-\int_0^1 \int_0^s\dot{E}_s(q(s))dtds\\
& &-\mathfrak{E}(q)-\frac{||\dot{q}||^2}{2}\\
&=&\mathfrak{E}^1(q)-\int_0^1 \frac{1}{|q(t)|}dt-\int_0^1 s\dot{E}_s(q(s))ds\\
&=&-\int_0^1 \frac{1}{|q(t)|}dt
\end{eqnarray*}	
which implies that $\mu=-1$. Hence in view of (\ref{mode}) we conclude that $q$ is a solution of (\ref{ode2}).

\section{Collisions}\label{coll}

In this section we suppose that $z \in \mathfrak{L^*Z}_0$ is a solution of (\ref{delay}). We do not assume anymore that $z \in \mathfrak{LZ}_0$. This means
the loop $z$ might pass sometimes through the origin which corresponds to collisions. We denote by
$$\mathfrak{C}_z:=\big\{\tau \in S^1: z(\tau)=0\big\}$$
the set of collision times. We first show that
\begin{equation}\label{finite}
|\mathfrak{C}_z|<\infty,
\end{equation}
i.e., the set of collision times is finite. For that purpose we first define $a_z \in \mathfrak{L}\mathbb{R}=C^\infty(S^1,\mathbb{R})$ for $\tau \in S^1$ by
\begin{eqnarray*}
a_z(\tau)&=&\frac{1}{2||z||^4}\bigg(2||z||^2||z'||^2+\mathcal{E}(z)+\mathcal{E}^1(z)-E_{t_z(\tau)}(z^2(\tau))\bigg)\\
& &-\frac{1}{2||z||^6}\bigg(\int_\tau^1 \dot{E}_{t_z(\sigma)}(z^2(\sigma))|z(\sigma)|^2d\sigma+1\bigg).
\end{eqnarray*}
We further define $b_z \in \mathfrak{L}\mathbb{C}$ for $\tau \in S^1$ by
$$b_z(\tau)=-\frac{|z(\tau)|^2}{2||z||^4}\nabla E_{t_z(\tau)}(z^2(\tau)),$$
and $c_z \in \mathfrak{L}(i\mathbb{R})$ by
$$c_z(\tau)=-\frac{i|z(\tau)|^2B(Z^2(\tau))}{||z||^2}.$$
Since $z$ solves (\ref{delay}) it follows that $z$ is as well a solution of the following second order linear homogeneous ODE
\begin{equation}\label{delayode}
z''(\tau)=a_z(\tau)z(\tau)+b_z(\tau)\bar{z}(\tau)+
c_z(\tau)z'(\tau).
\end{equation}
In order to prove (\ref{finite}) we show
\begin{equation}\label{nonvan}
z'(\tau_0) \neq 0, \quad \tau_0 \in \mathfrak{C}_z.
\end{equation}
Indeed, otherwise $z$ would be a solution of the second order ODE (\ref{delayode}) to the initial condition
$z(\tau_0)=z'(\tau_0)=0$ and therefore $z(\tau)=0$ for every $\tau \in S^1$. But this implies that $||z||=0$ contradicting the assumption that
$z \in \mathfrak{L^*Z}_0$. This shows (\ref{nonvan}). In particular $\mathfrak[{C}_z$ is a discrete subset of the circle and therefore
finite, so that (\ref{finite}) is proved as well. 
\\ \\
We define $t_z \colon S^1 \to S^1$ as in the regular case by
$$t_z(\tau)=\frac{\int_0^\tau|z(\sigma)|^2d\sigma}{||z||^2}, \quad \tau \in S^1.$$
Since the derivative of this function is given by
$$t'_z(\tau)=\frac{|z(\tau)|^2}{||z||^2}$$
we see that the set of collision times precisely coincides with the critical points of 
$t_z$, i.e.
$$\mathfrak{C}_z=\mathrm{crit}(t_z).$$
We next show that these critical points are inflection points, more precisely if $\tau_0 \in 
\mathfrak{C}_z$, than
\begin{equation}\label{inflection}
t_z'(\tau_0)=t_z''(\tau_0)=0, \quad t_z'''(\tau_0)>0.
\end{equation}
To see this we compute the second and third derivative of $t_z$ to be
$$t_z''(\tau)=\frac{2\langle z(\tau),z'(\tau)\rangle}{||z||^2}, \quad
t_z'''(\tau)=\frac{2\big(|z'(\tau)|^2+\langle z(\tau),z''(\tau)\rangle\big)}{||z||^2}, \qquad
\tau \in S^1.$$
If now $\tau_0 \in \mathfrak{C}_z$ we obtain in view of $z(\tau_0)=0$ and (\ref{nonvan}) that
$$t_z''(\tau_0)=0, \quad t_z'''(\tau_0)=\frac{2|z'(\tau_0)|^2}{||z||^2}>0$$
which proves (\ref{inflection}). In particular, $t_z$ is strictly monotone increasing and therefore gives rise to a homeomorphism from the circle to itself.
Therefore as in the regular case it has an inverse $\tau_z=t_z^{-1}$ which is still a homeomorphism of the circle. In contrast to the regular case
$\tau_z$ is not differentiable everywhere. However, on $S^1 \setminus t_z(\mathfrak{C}_z)$ it is still smoothly differentiable with derivative given by (\ref{dertau}). If we therefore define for $t \in S^1$ the function $q \colon S^1 \to \mathfrak{Q}_0$ by $q(t)=z^2(\tau_z(t))$, then $q$ is continuous and
smoothly differentiable on $S^1 \setminus t_z(\mathfrak{C}_z)$. 
\\ \\
In the following let us assume that $\mathfrak{C}_z \neq \emptyset$, so that $n:=|\mathfrak{C}_z|$ is a positive integer. Hence we can decompose the
complement of the set of collision times into connected components
$$S^1\setminus \mathfrak{C}_z=\bigcup_{j=1}^n \mathcal{I}_j$$
where each connected component
$$\mathcal{I}_j=(\tau_j^-,\tau_j^+)$$
is an open interval and
\begin{equation}\label{inter}
\tau_j^+=\tau_{j+1}^-, \quad 1 \leq j \leq n-1, \qquad \tau_n^+=\tau_1^-.
\end{equation}
For $1 \leq j \leq n$ we abbreviate
$$t_j^\pm:=t_z(\tau_j^{\pm})$$
and we set
$$\mathfrak{I}_j:=t_z(\mathcal{I}_j)=(t_j^-,t_j^+) \subset S^1.$$
The same argument as in the regular case shows that the restriction of $q$ to $\mathfrak{I}_j$ is a solution of (\ref{dotdot}) and there exists
$\mu_j \in \mathbb{R}$ such that 
\begin{equation}\label{muj}
\ddot{q}(t)=\frac{\mu_jq(t)}{|q(t)|^3}-B(q(t))i\dot{q}(t)-\nabla E_t(q(t)), \quad t \in \mathfrak{I}_j.
\end{equation}
To see this, we note that for $t \in \mathfrak{I}_j$ the second order delay equation (\ref{delayq}) still makes sense. In fact, since
the electric potential $E_t$ is also smooth at the origin the functionals $\mathfrak{E}$ and $\mathfrak{E}^1$ continuously extend by the same formula to $q$. Moreover, although $q$ can have zeros at which its derivative explodes in view of (\ref{kin}) and (\ref{pot}) the mean values on the righthand side of
(\ref{delayq}) are still finite. 
\\ \\
Combining (\ref{muj}) with (\ref{delayq}) we obtain for $\mu_j$ for any $t \in \mathfrak{I}_j$ the following equation
\begin{eqnarray}\label{muj2}
\mu_j&=&|q(t)|\bar{q}(t)\Big(\ddot{q}(t)+B(q(t))i\dot{q}(t)+\nabla E_t(q(t))\Big)\\ \nonumber
&=&|q(t)|\bigg(\mathfrak{E}(q)+\mathfrak{E}^1(q)+\frac{||\dot{q}||^2}{2}
-\int_0^1 \frac{1}{|q(s)|}ds-\int_t^1\dot{E}_s(q(s))ds-E_t(q(t))\bigg)\\ \nonumber
& &-\frac{|q(t)|\cdot|\dot{q}(t)|^2}{2}.
\end{eqnarray} 
Note that 
$$\lim_{t \to t_j^\pm}|q(t)|\bigg(\mathfrak{E}(q)+\mathfrak{E}^1(q)+\frac{||\dot{q}||^2}{2}
-\int_0^1 \frac{1}{|q(s)|}ds-\int_t^1\dot{E}_s(q(s))ds-E_t(q(t))\bigg)=0$$
since all terms in the bracket remain finite as $t$ goes to $t_j^\pm$. Using (\ref{der1}) and the definition of $q$ we have
$$\frac{|q(t)|\cdot|\dot{q}(t)|^2}{2}=\frac{4|z^2(\tau_z(t))|^2||z||^4||z'(\tau_z(t))||^2}{2|z^2(\tau_z(t))|^2}=2||z||^4||z'(\tau_z(t))||^2.$$
Therefore if we take in (\ref{muj2}) the limit $t \to t_j^\pm$ we obtain
\begin{equation}\label{muj3}
\mu_j=-2||z||^4|z'(\tau_j^+)|^2=-2||z||^4|z'(\tau_j^-)|^2, \quad 1 \leq j \leq n.
\end{equation}
In view of (\ref{inter}) we conclude that all $\mu_j$ are equal, i.e., there exists $\mu$ independent of $j$ such that
$$\mu_j=\mu, \quad 1 \leq j \leq n.$$
Hence we obtain from (\ref{muj}) that
$$\ddot{q}(t)=\frac{\mu q(t)}{|q(t)|^3}-B(q(t))i\dot{q}(t)-\nabla E_t(q(t)), \quad t \in S^1 \setminus t_z(\mathfrak{C}_z).$$
As in the regular case we conclude from this equation that $\mu=-1$.

\section{Correspondence of solutions}\label{correspondence}

The complex squaring map $\varsigma \colon \mathfrak{Z} \to \mathfrak{Q}$ is a double cover. This has several consequences. First we observe that the functional
$\mathcal{F} \colon \mathfrak{L^*Z}_0 \to \mathbb{R}$ is invariant under the involution
$$I \mathfrak{L^*Z}_0 \to \mathfrak{L^*Z}, \quad z \mapsto -z,$$
i.e. we have
$$\mathcal{F} \circ I=\mathfrak{F}.$$
This means that if $z$ is a critical point of $\mathcal{F}$, then $-z$ is as well a critical point. In other words if $z \in \mathfrak{L^*Z}_0$ is a solution
of the second order delay equation (\ref{delay}), than the same is true for $-z$, as one can see as well directly by looking at (\ref{delay}). Both $z$ and
$-z$ correspond to the same (collisional) solution of (\ref{ode2}). 
\\ \\
The loop space $\mathfrak{LQ}$ decomposes as
$$\mathfrak{LQ}=\bigcup_{n \in \mathbb{Z}}\mathfrak{L}_n\mathfrak{Q}$$
where $\mathfrak{L}_n\mathfrak{Q}$ consists of loops in $\mathfrak{LQ}$ whose winding number around the origin is $n$. 
The same decomposition we have as well for the loop space of $\mathfrak{Z}$, i.e.,
$$\mathfrak{LZ}=\bigcup_{n \in \mathbb{Z}}\mathfrak{L}_n\mathfrak{Z}.$$
The complex squaring map doubles the winding number so that $\Sigma$ maps $\mathfrak{L}_n\mathfrak{Z}$ to $\mathfrak{L}_{2n}\mathfrak{Q}$ and $\Sigma$ induces a diffeomorphism
$$\Sigma_n \colon \mathfrak{L}_n\mathfrak{Z}/I \to \mathfrak{L}_{2n}\mathfrak{Q}, \quad n \in \mathbb{Z}.$$
To take into account as well loops in $\mathfrak{LQ}$ of odd winding number we consider the twisted loop space
$$\mathfrak{L}_-\mathfrak{Z}:=\big\{z \in C^\infty(\mathbb{R},\mathfrak{Z}):z(t+1)=-z(t),\,\, t\in \mathbb{R}\big\}.$$
Note that the twisted loop space $\mathfrak{L}_-$ is still invariant under the involution $I$ which maps $z$ to $-z$. We further have the doubling map
$$\Delta \colon \mathfrak{L}_-\mathfrak{Z} \to \mathfrak{LZ}$$
which for $z \in \mathfrak{L}_-\mathfrak{Z}$ is given by
$$\Delta(z)(t)=z(2t), \quad t \in \mathbb{R}.$$
Indeed, note that
$$\Delta(z)(t+1)=z(2t+2)=-z(2t+1)=z(2t)=\Delta(z)(t)$$
so that $\Delta(z)$ is periodic with period $1$. Therefore, for an odd integer $n$ we can set
$$\mathfrak{L}_{\frac{n}{2}}\mathfrak{Z}=\big\{z \in \mathfrak{L}_-\mathfrak{Z}: \Delta(z) \in \mathfrak{L}_n\mathfrak{Z}\big\}$$
so that we have
$$\mathfrak{L}_- \mathfrak{Z}=\bigcup_{n \in \mathbb{Z}}\mathfrak{L}_{\frac{1}{2}+n}\mathfrak{Z}.$$
In particular, we can think of twisted loops as loops having winding number in $\mathbb{Z}+\tfrac{1}{2}$. 
We define
$$\Sigma \colon \mathfrak{L}_-\mathfrak{Z} \to \mathfrak{LQ}$$
by the same formula as in the untwisted case. Then we have induced diffeomorphisms
$$\Sigma_n \colon \mathfrak{L}_n\mathfrak{Z}/I \to \mathfrak{L}_{2n}\mathfrak{Q}, \quad n \in \frac{1}{2}\mathbb{Z}$$
for integer and half-integer winding numbers. 
\\ \\
After having extended $\Sigma$ as well to twisted loops we can pullback the functional $\mathfrak{F}\colon \mathfrak{LQ} \to  \mathbb{R}$ to
$$\Sigma^* \mathfrak{F} \colon \mathfrak{LZ} \cup \mathfrak{L}_-\mathfrak{Z} \to \mathbb{R}.$$
The formula (\ref{pullback}) continuous to hold for twisted loops. After extension $\Sigma$ gives rise to a one-to-one correspondence between critical points
$$\mathrm{crit} \Sigma^* \mathfrak{F}/I \cong \mathrm{crit}\mathfrak{F}$$
where critical points of $\mathfrak{F}$ of winding number $n$ correspond to critical points of $\Sigma^*\mathfrak{F}$ of winding number
$\tfrac{n}{2}$. 
\\ \\
In the same way as we have blown up the space of untwisted loops we can as well blow up the space of twisted loops. Namely we set 
$$\mathfrak{L}^*_-\mathfrak{Z}_0=\big\{z \in \mathfrak{L}_-\mathfrak{Z}_0: ||z||>0\big\}.$$
The pullback of $\mathfrak{F}$ under $\Sigma$ extends as well smoothly to the blow up of twisted loops by the same formula (\ref{pullback}), so that we can consider the functional $\mathcal{E}$ as a functional
$$\mathcal{E} \colon \mathfrak{L}^*\mathfrak{Z}_0 \cup \mathfrak{L}^*_-\mathfrak{Z}_0 \to \mathbb{R}$$
which is defined on the blow up of twisted and untwisted loops. In the twisted and untwisted case critical point are solutions of the second order delay equation (\ref{delay}). By the same argument as in the untwisted case also in the twisted case $\Sigma$ maps critical points of
$\mathcal{E}$ to collisional periodic solutions of (\ref{ode2}). Although for collisional solutions of (\ref{ode2}) we cannot associate anymore a winding number we are still able to associate to them a \emph{winding number modulo 2}. Namely we say that a collisional solution of (\ref{ode2}) has
even winding number if it corresponds under $\Sigma$ to a critical point of $\mathcal{E}$ on
$\mathfrak{L}^*\mathfrak{Z}_0$ and it has odd winding number if it corresponds to a critical point on $\mathfrak{L}^*_-\mathfrak{Z}_0$. Hence we can summarize the findings of this note in the following theorem.
\begin{thm}\label{main}
There is a one-to-one correspondence between $\mathrm{crit}(\mathcal{E})/I$ and collisional
periodic solutions of the ODE (\ref{ode2}). Untwisted critical points correspond to collisional solutions of even winding number and twisted critical points to collisional solutions of odd winding number. Moreover, critical points of $\mathcal{E}$ are solutions of the second order delay equation (\ref{delay}).
\end{thm}

\end{document}